\documentclass[12pt,a4paper]{article}
\usepackage{fleqn,amssymb,amsmath,amscd,epsfig,version,theorem} %
\usepackage{hyperref}
\newcommand{\Section}[1]{\section{#1}\setcounter{equation}{0}}
\numberwithin{equation}{section}
\numberwithin{figure}{section}
\theoremstyle{change}
\newtheorem{theorem}{Theorem} [section]
\newtheorem {lemma}[theorem]{Lemma}

\newtheorem {defi}[theorem]{Definition}
{\theorembodyfont{\normalfont}
{\theorembodyfont{\normalfont}\newtheorem {example}[theorem]{Example}
{\theorembodyfont{\normalfont}
{\theorembodyfont{\normalfont}\newtheorem {remark}[theorem]{Remark}
{\theorembodyfont{\normalfont}\newtheorem {remarks}[theorem]{Remarks}
\newcommand{\beq}{\begin{equation}}
\newcommand{\eeq}{\end{equation}}
\newcommand{\Leq}[1]{\label{#1}\end{equation}}
\newcommand{\beqn}{\begin{eqnarray}}
\newcommand{\eeqn}{\end{eqnarray}}
\newcommand{\beqno}{\begin{eqnarray*}}
\newcommand{\eeqno}{\end{eqnarray*}}

\renewcommand {\l}{\left}
\newcommand {\ri}{\right}
     
\newcommand {\LA}{\left\langle}
\newcommand {\RA}{\right\rangle}
\newcommand {\pa}{\partial}

\newcommand {\eh}{{\textstyle \frac{1}{2}}}

\newcommand {\sign}{{\rm sign}}

\newcommand {\bJ}{{\mathbb J}}

\newcommand {\bR}{{\mathbb R}}

\newcommand {\bZ}{{\mathbb Z}}

\newcommand{\idty}{{\rm 1\mskip-4mu l}} 
\newcommand{\ov}{\overline}
\newcommand{\bem}{\l(\! \begin{array}}
\newcommand{\eem}{\end{array}\!\ri)}
\newcommand{\bsm}{\left(\begin{smallmatrix}} 
\newcommand{\esm}{\end{smallmatrix}\right)}  

\newcommand{\qmbox}[1]{\quad\mbox{#1}\quad}
\renewcommand {\max}{{{\rm max}}}

\def\XXint#1#2#3{{\setbox0=\hbox{$#1{#2#3}{\int}$}
\vcenter{\hbox{$#2#3$}}\kern-.5\wd0}}

\usepackage{color}

\begin{document}
\title{Symbolic Dynamics of Magnetic Bumps}
\author{Andreas Knauf\thanks{Department Mathematik, 
Universit\"{a}t Erlangen-N\"{u}rnberg,
Cauerstr. 11, D--91058 Erlangen, Germany. 
e-mail: \texttt{knauf@math.fau.de}} , Marcello Seri\thanks{Department of Mathematics and Statistics
University of Reading
Whiteknights, PO Box 220, Reading RG6 6AX (UK)
e-mail: \texttt{m.seri@reading.ac.uk}}}
\date{\today}
\maketitle
\begin{abstract}
For $n$ convex magnetic bumps in the plane, whose boundary has a curvature somewhat smaller
than the absolute value of the constant magnetic field inside the bump, we 
construct a complete symbolic dynamics of a classical particle moving with speed one.

\end{abstract}
%
%
%
\Section{Introduction and notation}
The subject of chaotic scattering is mostly about scattering by obstacles and by
potential bumps, see \cite{Ga,Sm}.
Here we consider the motion of a classical particle in the plane under the influence
of a magnetic potential. In the case of motion in the plane, this has the form
\[B:= \sum_{\ell\in A} B_\ell\qmbox{with}B_\ell:=b_\ell\, \idty_{C_\ell}:\bR^2\to \bR,\] 
with the alphabet $A:=\{1,\ldots,n\}$, for mutually
disjoint, convex and compact domains $C_\ell\subseteq \bR^2$ with $C^2$ boundaries 
$\partial C_\ell$, and field strengths $b_\ell\in\bR\setminus \{0\}$.
The Hamiltonian system $(P,\omega,H)$ with phase space $P:= T\bR^2$, 
\[H:P\to\bR\qmbox{,} H(q,v)=\eh \|v\|^2\]
and (discontinuous) symplectic form 
\[\omega:= dq_1\wedge dv_1+dq_2\wedge dv_2 + B\,dq_1\wedge dq_2\]
give rise to the Hamiltonian vector field $X:P\to TP$, ${\mathbf i}_X\omega =dH$.
The corresponding differential equation is
\[\dot{q} = v \qmbox{,}  \dot{v} = B(q)\,\bJ\, v, \qquad \mbox{ with } \bJ := \bsm0&-1\\1&0\esm.\]
The solution of the initial value problem for magnetic field $b\in\bR\setminus\{0\}$ thus equals
\[q(t) = q_0 + \frac{1}{b} 
\bsm \sin(bt)&\cos(bt)-1\\
1-\cos(bt)&\sin(bt)\esm v_0\qmbox{,}v(t) = \bsm \cos (bt)&-\sin (bt)\\ \sin (bt)&\cos (bt)\esm v_0.\]
So  the particle with energy $H=E$
has speed $\sqrt{2E}$. Outside the bumps it moves on a straight line, whereas
inside the $\ell^{\rm th}$ bump it moves on a segment of a circle of Larmor radius $\sqrt{2E}/|b_\ell|$.
The sense of rotation is positive (counter-clockwise), if $b_\ell>0$.
The $t$--invariant center of the circle is 
\beq
c_\ell(v,q):=q+b_\ell^{-1}\bJ v.
\Leq{Larmor:center}

Without loss of generality we fix the energy to equal $E:=1/2$, so that we get unit speed
on the energy surface 
\[\hat\Sigma:=H^{-1}(E)\qmbox{with projection} \pi:\hat\Sigma\to\bR^2,\ (q,v)\mapsto q .\]
\section{Single bumps}
\begin{figure}[htbp]
\begin{center}
\centerline{\includegraphics[width=110mm,natwidth=529,natheight=464]{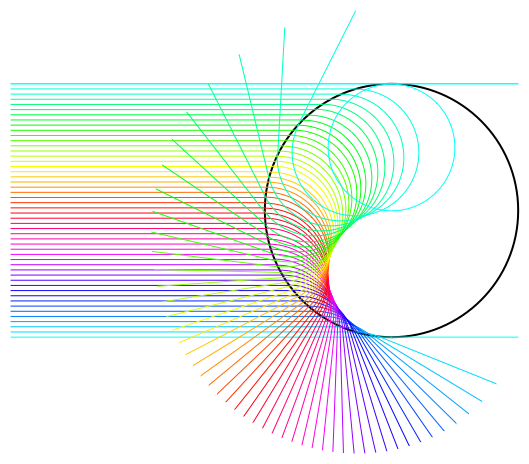}}
\caption{Dynamics for a disk $C$ and a strong field (with $b=-2r$)}
\label{fig:strong}
\end{center}
\end{figure}

\begin{figure}[htbp]
\begin{center}
\centerline{\includegraphics[width=120mm,natwidth=622,natheight=358]{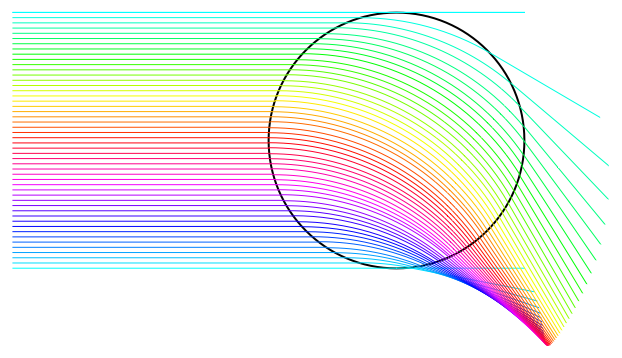}}
\caption{Dynamics for a disk $C$ and a weak field (with $b=-r/2$)}
\label{fig:weak}
\end{center}
\end{figure}

Consider for $\ell\in A$ the set $\tilde{C}_\ell$ of
points  $c\in C_\ell$ whose minimal distance to the boundary $\partial C_\ell$ 
is larger than $|b_\ell|^{-1}$. This is a compact convex set, possibly void.\footnote{If 
$\tilde{C}_\ell$ has positive area, its boundary is Lipschitz. However, e.g.\ for a stadion $C_\ell$
composed of two half-disks of radius $1/|b_\ell|$ and a rectangle,
$\tilde{C}_\ell$ is a  line segment.}
We remove from $\hat\Sigma$ the sets
$\Sigma_\ell:=\{(v,q)\in \hat\Sigma\mid c_\ell(v,q)\in \tilde{C}_\ell\}$. 
These compact sets are invariant
under the $2\pi/|b_\ell|$--periodic Larmor flow for magnetic field $b_\ell$. 
At points $x\in \partial \Sigma_\ell$ corresponding to Larmor circles in the configuration plane that
touch $\partial{C}_\ell$ in only one point, the boundary
$\partial \Sigma_\ell$ is $C^2$. 
In particular this is the case if the curvature $\kappa_\ell\in C(\partial  C_\ell,[0,\infty))$
of $\partial  C_\ell$ is strictly smaller than $|b_\ell|$. 
\begin{defi}\quad\\
$\bullet$
We set $\underline{\kappa}_\ell:=\min_{x\in \partial  C_\ell}\,\kappa_\ell(x)$ and
$\overline{\kappa}_\ell:=\max_{x\in \partial  C_\ell}\,\kappa_\ell(x)$ 
\ (so $0\le \underline{\kappa}_\ell \le\overline{\kappa}_\ell > 0$).\\
$\bullet$
If $|b_\ell| < \underline{\kappa}_\ell$ for all $\ell\in A$, the system has {\bf weak fields}.\\
$\bullet$
If $|b_\ell| > \overline{\kappa}_\ell$ for all $\ell\in A$, the system has {\bf strong fields}.\\
$\bullet$
For $n\ge 2$ we set 
$d_\ell:=\min_{m\in A\setminus\{\ell\}}\,\min_{x\in C_\ell, y\in C_m}\|x-y\|$
\ (so $d_\ell>0$).\\
$\bullet$
If $|b_\ell| >  1/(d_\ell\, \alpha_{\min})+2\overline{\kappa}_\ell $ for all $\ell\in A$, the system has 
{\bf very strong fields}.\footnote{Here $\alpha_{\min}$ is a constant, related to the arrangement
of the bumps in the plane, that will be discussed in more details in Section \ref{sect:cone}.}
\end{defi}
\begin{example}[Single disks]\quad\\  
For $n=1$ and a disk $C\equiv C_1$ of radius $r>0$ the curvature of $\partial C$ 
is constant and equals $1/r$. So the field $B=b\, \idty_{C}$ is weak iff $|b|<1/r$ and strong iff $|b|>1/r$.
Compare Figures \ref{fig:strong} and \ref{fig:weak} for dynamics near such a bump.
\hfill $\Diamond$
\end{example}
We will mainly consider the (very) strong fields cases.

\begin{figure}[htbp]
\begin{center}
\centerline{\includegraphics[height=60mm,natwidth=529,natheight=464]{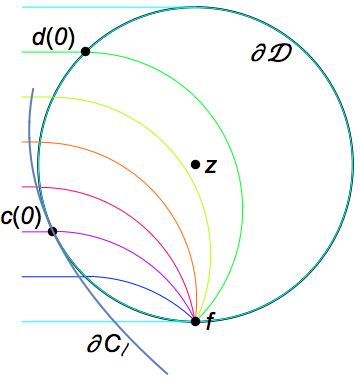}
\hfill\includegraphics[height=60mm,natwidth=360,natheight=352]{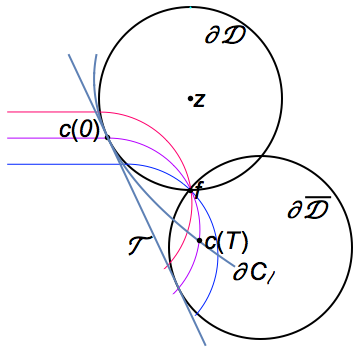}}
\caption{Proof of Lemma \ref{lem:focus}}
\label{fig:lemma}
\end{center}
\end{figure}

The following lemma relates strong fields to some property of a single bump.
Our convention for the 'Jacobi equation'\,\footnote{This 
is an abuse of language, since Hamilton's equation is not geodesic.} along the trajectory
is to use the orthonormal basis $e_1(t),\,e_2(t)$ of $T_{c(t)}\bR^2$, given by 
$e_1(t):=v(t)$ and $e_2(t):=\bJ e_1(t)$. 
Then writing a vector field along $c$
in the form $t\mapsto I(t)e_1(t)+J(t)e_2(t)$, we obtain the linearized flow with
\[J(t)=-\sin (bt)I(0)+\cos (bt)J(0)- (1-\cos (bt))\dot{I}(0)/b+\sin(bt)\dot{J}(0)/b\]
and
\[\dot{J}(t)=-b\cos (bt)I(0)-b\sin (bt)J(0)- \sin (bt)\dot{I}(0)+\cos(bt)\dot{J}(0).\]
\begin{lemma}  \label{lem:focus}\quad\\
In the case of a strong field, along a trajectory $c:I\to \bR^2$ 
entering the bump at $c(0)\in \partial  C_\ell$ and leaving it at time $T>0$
($c(T)\in \partial  C_\ell$),
the parallel incoming Jacobi field\,\footnote{We assume w.l.o.g.\ that the component 
parallel to the direction $\dot{c}(0)$
vanishes.} $\bigl(J(0),\dot{J}(0)\bigr)=(1,0)$
has outgoing data  $J(T)<0$, $\dot{J}(T)<0$.
\end{lemma}
\textbf{Proof:}
$\bullet$
Consider the unique disk $\cal D$ of radius $r:=1/ |b_\ell|$, whose boundary is tangent to 
$\partial  C_\ell$ at $c(0)$
and which is entered by the trajectory at time $t=0$. As the curvature $|b_\ell|$ of its boundary
is strictly larger than the curvature of $\partial  C_\ell$, we have ${\cal D}\subseteq  C_\ell$,
and $c(0)$ is the only intersection between $\partial {\cal D}$ and $\partial C_\ell$.
We denote the center of ${\cal D}$ by $z\in\bR^2$, see Figure \ref{fig:lemma}, left.\\
$\bullet$
We claim that all trajectories $d:I\to \bR^2$ entering ${\cal D}$ at time $t=0$ with velocity 
$\dot {d}(0) = \dot {c}(0)$ leave  ${\cal D}$ at the same point $f:=z+\bJ\dot {c}(0)/b_\ell$ :\\
Clearly $d(0)$ is a solution $x$ of the equations
\[\|x-z\|^2=r^2\qmbox{and}\|x-d(0)-\bJ\dot {c}(0)/b_\ell \|^2=r^2\]
of $\partial {\cal D}$ and the Larmor circle \eqref{Larmor:center}.
But the second intersection of these two circles equals $f$.\\
$\bullet$
The union of the above orbit family is a Lagrangian manifold, and the vector in 
$T_{x} \Sigma$, $x:=(c(0),\dot{c}(0))$ given Jacobi field data
$\bigl(J(0),\dot{J}(0)\bigr)$, is tangent to it. So the Jacobi field along $c$ turns
vertical at $f=c(t_f)$, that is, $J(t_f)=0$ and $\dot{J}(t_f)<0$. 
But $f$, being a boundary point of ${\cal D}$ different from $c(0)$,
is an interior point of $C_\ell$, so that $t_f\in (0,T)$.\\
$\bullet$
Consider the line ${\cal T}$ in $\bR^2$ that is tangent  to $\partial  C_\ell$ at $c(0)$.
If we denote by $\ov{\cal D}$ the image of ${\cal D}$ under the reflection by the line
through $f$, perpendicular to ${\cal T}$ (see Figure \ref{fig:lemma}, right), then 
the Larmor circle through $c(0)$ and $f$ intersects $T$ for the second time at
the unique point in $T\cap \pa \ov{\cal D}$.
Again by reflection symmetry, at that point the Jacobi field along $c$ is parallel.\\
$\bullet$
As $C_\ell$ is strictly convex, it lies entirely on one side of ${\cal T}$. So by a comparison argument
(see again Figure \ref{fig:lemma}, right),
$J(T)<0$.
\hfill $\Box$
\begin{remark} 
If $C_\ell$ is a disk, in the strong field case we conclude from Lemma  \ref{lem:focus},
that for parallel incoming trajectories the envelope of the solution curves  is
a half-circle of radius $r_\ell-1/|b_\ell|$, see lower part of Figure \ref{fig:strong}.
\hfill $\Diamond$
\end{remark}
\section {The degree for weak and strong fields}
Scattering by a bump is qualitatively different for weak and for strong force fields.
To see this, consider oriented lines in configuration space $\bR^2$.
The set of these lines is naturally isomorphic to the cylinder $TS^1\cong S^1\times \bR$.

We assume that there is just one bump (so $n=1$), and consequently omit subindices.
There is a problem in defining the flow on $\Sigma$ in the case of {\em glancing} trajectories,
that is, trajectories tangent to $\pa C$ (consider the uppermost trajectory in Figure \ref{fig:strong}). 
Either we continue these incoming rays 
by just extending them to a straight
line, or we extend them by a segment of the Larmor circle in the bump. This is a complete
circle for a strong field $b$. So in that case we could attach the outgoing ray with the incoming
direction.

None of these prescriptions leads to a continuous {\em flow}.
However, it is important to notice that  either prescription leads to a continuous {\em map} called
\[S:TS^1 \to TS^1,\]
sending incoming to outgoing oriented lines.\,\footnote{For fields that are neither weak nor strong,
it can happen that $S$ cannot be defined continuously.}

As shown in \cite{Kn}, see also \cite{KK}, such a map defines a topological index 
${\rm deg}(C,b)\in\bZ$.  Using the bundle projection
$\pi: TS^1\to S^1$, it can be defined as follows.
A family of incoming rays with the same direction is parameterized by the value of 
their angular momentum
$L:\Sigma\to \bR,\ (q,v)\mapsto \LA \bJ q,v\RA$. This is obviously flow-invariant outside
the bump. $S$ maps initial to final pairs, consisting of directions in $S^1$ and 
angular momenta in $\bR$.

Given the initial direction $\varphi \in S^1$, the continuous map 
\[T_\varphi S^1\cong \bR\to S^1\qmbox{,}\ell\to \pi\circ S(\ell, \varphi)\] 
can be uniquely completed
(using the Alexandrov compactification $\bR\cup\{\infty\}\cong S^1$ of $\bR$) to
a continuous map $S^1\to S^1$. The degree of that map is independent of $\varphi\in S^1$ and 
is called the {\em scattering degree}.
\begin{lemma}  \label{lem:deg}
For a single bump $C$, the scattering degree ${\rm deg}(C,b)$ equals
\begin{enumerate}
\item[$\bullet$] 
zero, for weak fields $b$,
\item[$\bullet$] 
$\sign(b)$, for strong fields $b$.
\end{enumerate}
\end{lemma}
\textbf{Proof:}
The total curvature of a trajectory equals the signed angular length of the segment of its
Larmor circle. This is obviously zero if the trajectory does not intersect $C$.
It is bounded away from $2\pi$ in absolute value for weak fields.\\
For strong fields $b>0$ it equals $+2\pi$ exactly for the glancing trajectories intersecting $\pa C$
with $C$ on their left hand side. It then decreases to zero when the angular momentum value
of the incoming rays it decreased, until one arrives at a glancing trajectory with $C$ on its right
hand side.\\ 
For strong fields $b<0$ the sides are interchanged.
\hfill $\Box$
\section {A strictly invariant cone field}\label{sect:cone}
\noindent
{\bf Assumption:}\\ 
We assume that no three convex sets $C_\ell\subseteq \bR^2$ lie on a line.
For $n\ge3$
\[\alpha_{\min}:=\!\! \min_{k\neq \ell\neq m\neq k\in  A}\!\!
\min \Big\{ \arccos \Big( \! \LA {\textstyle \frac{y-x}{\|y-x\|}\, ,\, \frac{y-z}{\|y-z\|}} \! \RA \Big)\mid
x\in C_k, y\in C_\ell, z\in C_m\Big\},\]
and for $n=2$  we set $\alpha_{\min}:=\pi/3$, say.
As we assumed that no three bumps are intersected by a line, the angle $\alpha_{\min}\in (0,\pi/3]$ 
is positive.
For a segment of a solution curve that intersects $C_k$, $C_\ell$ and $C_m$ in succession,
the total curvature inside $C_\ell$ is bounded below by $\alpha_{\min}$. 
So the length of that sub-segment is at least $\alpha_{\min}/|b_m|$.

We now consider in the very strong force regime the set 
$\Lambda\subseteq \Sigma$ of initial  conditions
belonging to bounded orbits. We denote by
\[{\cal N}(x)\in T_x\bR^2\qquad (\ell\in A, x\in \pa C_\ell) \]
the unit outward normal vectors of $C_\ell$.\\ 
${\cal H}^\pm:=\bigcup_{\ell\in A} {\cal H}_\ell \subseteq \Sigma$ 
is the disjoint union of the Poincar\'e surfaces 
\[ {\cal H}_\ell^\pm := \{ (v,q)\in \Sigma\mid q\in \partial C_\ell, \ \pm \LA v , {\cal N}(q)\RA>0\}.\]
For arbitrary fields, the flow induces a Poincar\'e
map ${\cal P}^{(i)}: {\cal H}^-\to {\cal H}^+$ internal to the bumps, which is a diffeomorphism.
On the other hand, there are maximal open subsets $V^\pm\subseteq {\cal H}^\pm$
so that the flow gives rise to a diffeomorphism 
${\cal P}^{(e)}:V^+\to V^-$, external to the bumps. 
Setting $U^-:= ({\cal P}^{(i)})^{-1}(V^+)$, we obtain a diffeo
\[{\cal P}:={\cal P}^{(e)}\circ {\cal P}^{(i)}: U^-\to V^-.\]

\begin{figure}[htbp]
\begin{center}
\centerline{\includegraphics[height=60mm,natwidth=360,natheight=289]{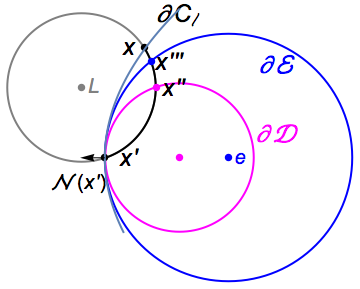}}
\caption{Proof of Lemma \ref{lem:conefield}}
\label{fig:lemma:cone}
\end{center}
\end{figure}

\begin{lemma}  \label{lem:conefield}
Assume that the fields $b_\ell$ obey the very strong fields inequalities
\beq 
|b_\ell|  > 1/(d_\ell\, \alpha_{\min})+2 \ov{\kappa}_\ell \quad (\ell\in A). 
\Leq{ineq:curv}
Then the cone fields 
\[{\cal C}_\ell(x):= \{ \lambda^{(l)}e^{(l)}+ \lambda^{(u)}e^{(u)}\mid \lambda^{(l)},\lambda^{(u)}\in\bR, 
\lambda^{(l)}\cdot \lambda^{(u)}>0  \}\quad (x\in {\cal H}_\ell),\]
defined by the tangent vectors $e^{(l)},e^{(u)}\in T{\cal H}_\ell$,
\[e^{(l)}:= \bsm1\\0\esm\qmbox{,} e^{(u)}:=\bsm d_{\ell}\\ 1\esm\] 
are strictly invariant under the linearized Poincar\'e map $T{\cal P}$.
\end{lemma}
\textbf{Proof:}\\
$\bullet$
According to Lemma \ref{lem:focus}, the vector $e^{(l)}\in T_xU^-$ (with $x\in {\cal H}_\ell^-$) 
is mapped by $T_x{\cal P}$ 
to a vector $T_y{\cal P}(e^{(l)})\in T_{y}V^-$ (with $y:={\cal P}(x)\in {\cal H}_m^-$)
of the form $\bsm 1&d\\ 0&1\esm \bsm J\\ \dot{J}\esm$ with $J,\dot{J}<0$, $d$ being the
length of the trajectory segment between the points $\pi\circ{\cal P}^{(i)}(x)$  of exit from $C_\ell$
resp.\ $\pi\circ{\cal P}(x)$ of entrance in $C_m$. \\
This shows that $e^{(l)}$ is contained in the cone field ${\cal C}_\ell(x)$.\\
$\bullet$
We want to show a similar statement for the vector $e^{(u)}\in T_xU^-$, but we choose to work
backwards in time. So we consider the vector $e^{(u)}=\bsm d_{m}\\ 1\esm\in T_yU^-$ and show that its preimage
$T_y({\cal P})^{-1}(e^{(l)})$ is {\em not} contained in ${\cal C}_\ell(x)$. 
First of all, $T_y({\cal P}^{(e)})^{-1}(e^{(l)})=\bsm d_m-d\\ 1\esm$ has non-positive first entry, since $T_y({\cal P}^{(e)})^{-1}=
\bsm 1&-d\\ 0&1\esm$ and $d\ge d_m$. For comparison we consider the vector 
$\bsm 0\\ 1\esm\in T_{x'}$ with $x':={\cal P}^{(i)}(x)$ and show that its preimage 
$T_{x'}({\cal P}^{(i)})^{-1}\big(\bsm 0\\ 1\esm\big)$ is not contained in  ${\cal C}_\ell(x)$.\\
We choose the disk ${\cal E}\subseteq C_\ell$ of radius $1/\overline{\kappa}_\ell$, whose boundary 
is tangent to $\pa C_\ell$ at $x'$. Likewise, ${\cal D}\subseteq {\cal E}$ is the disk of radius $1/|b_\ell|$
whose boundary is tangent to $\pa C_\ell$ (and $\pa {\cal E}$) at $x'$; see Figure \ref{fig:lemma:cone}.\\
Like in the proof of Lemma \ref{lem:focus}, at the second intersection $x''$ of $\pa {\cal D}$ with the Larmor circle (gray in Figure \ref{fig:lemma:cone}), the family of Larmor solutions crossing at $x'$
has become parallel. We must show that by following the corresponding Jacobi field backwards
from $x''$ to $x$, it is not contained in the cone field ${\cal C}(x)$. 
Although we have no direct information about the length of the Larmor circle segment between
$x''$ and $x$, we know (by definition of $\ov{\kappa}_\ell $) 
that it is longer than the one between $x''$ and
the intersection $x'''$ with $\pa {\cal E}$.  

So we must bound the length of that arc from below.
The length of the arc between $x'$ and $x'''$ equals $|b_\ell| \alpha$, $\alpha$ being the angle 
between the normal ${\cal N}(x')$ and the center $L$ of the Larmor circle, seen from $x'$.

The Larmor angle between $x'$ and $x'''$ equals $2 \arctan[r \sin(\alpha)/(1 + r \cos(\alpha))]$, with
$r:=\ov{\kappa}_\ell / |b_\ell|$ being the ratio between the radii of $\pa {\cal D}$ and $\pa {\cal E}$.
This identity follows when considering the line through $L$ and the center $e$ of ${\cal E}$.
That line bisects the Larmor angle. By elementary trigonometry the angle between that line
and the one through $e$ and $x'$ equals  $\arctan[r \sin(\alpha)/(1 + r \cos(\alpha))]$.

Thus the length of the Larmor arc between $x''$ and $x'''$ equals 
\[f(\alpha):=\alpha-2 \arctan[r \sin(\alpha)/(1 + r \cos(\alpha))],\] 
multiplied by its radius $1/|b_\ell|$.
As $f(0)=0$, $f'(\alpha)=\frac{1-r^2}{r^2+2 r \cos (\alpha )+1}$ with $f'(0)=\frac{1-r}{1+r}>0$
and $f''(\alpha)=\frac{r \left(1-r^2\right) \sin (\alpha )}{\left(r^2+2 r \cos (\alpha )+1\right)^2}\ge0$, 
we can bound $f$ from below by
\beq
f(\alpha)\ge {\textstyle \frac{1-r}{1+r}}\;\alpha=  
{\textstyle \frac{|b_\ell|-\ov{\kappa}_\ell}{|b_\ell|+\ov{\kappa}_\ell}}\;\alpha\qquad (\alpha\in[0,\pi]).
\Leq{f:ineq}
As the Larmor arc belongs to a solution segment that intersects $C_k$, $C_\ell$ and $C_m$ in
consecution, the argument of $f$ is bounded below by $\alpha \ge \alpha_{\min}$.

We have to check, whether then the inequality
\[d_\ell^{-1}<\tan(|b_\ell|\,f(\alpha)) \] 
holds true. By \eqref{f:ineq} we check the stronger inequality
\[d_\ell^{-1}<\tan\l(|b_\ell|\,{\textstyle \frac{|b_\ell|-\ov{\kappa}_\ell}{|b_\ell|+\ov{\kappa}_\ell}}\;\alpha\ri). \]
In fact, we argue that even
\[d_\ell^{-1}<|b_\ell|\,{\textstyle \frac{|b_\ell|-\ov{\kappa}_\ell}{|b_\ell|+\ov{\kappa}_\ell}}\;\alpha. \]
The corresponding equation is quadratic in $x:=|b_\ell|$ and has the positive solution
$x=\eh(c+ \ov{\kappa}_\ell+\sqrt{c^2+6c \ov{\kappa}_\ell+\ov{\kappa}_\ell^2})\le c+2 \ov{\kappa}_\ell$
with $c:=1/(d_\ell\alpha_{\min})$.
\hfill $\Box$\\[2mm]
\section {Symbolic dynamics}
We equip for the (discretely topologized) alphabet $A=\{1,\ldots,n\}$ the sequence space
\[\Xi_A:=\{a:\bZ \to A\mid \forall k\in \bZ: a_{k+1}\neq a_k\}\]
with the product topology. Then the shift $\sigma:\Xi_A \to \Xi_A$, $\sigma(a)_k=a_{k+1}$
is known to be a homeomorphism.
\begin{theorem}  
Let the very strong fields assumption \eqref{ineq:curv} be valid for the $n\ge2$ bumps.
Then every bounded trajectory $c:\bR\to\Sigma$ intersects ${\cal H}$ infinitely  often.
The set $\Xi:=\Lambda \cap {\cal H}$ is homeomorphic to $\Xi_A$,
under the homeomorphism\,\footnote{identifying one-element sets with their elements}
\beq
\Phi: \Xi_A\to \Xi\qmbox{,} a\mapsto \bigcap_{k\in \bZ} \l\{ {\cal P}^k(x)\mid x\in {\cal H}_{a_{0}}: 
{\cal P}^k(x) \mbox { is defined and in }{\cal H}_{a_{k}}\ri\} 
\Leq{def:Phi}
to one-point sets. $\Phi$ is a conjugacy for the Poincar\'e map ${\cal P}$ and the shift map 
$\sigma$, i.e.\ 
\[\Phi\circ \sigma = {\cal P}\circ \Phi.\]
\end{theorem}
\textbf{Proof:}
$\bullet$
It is generally true that for a continuous, strictly invariant cone field there is at most one element
in each set $\Phi(a)$ of \eqref{def:Phi}.\\
$\bullet$
Then we employ the technique of \cite[Chapter 6]{Kn} 
to show that the sets defining $\Phi(a)$ are non-void. According to Lemma \ref{lem:deg}
we have degrees $\deg(C_\ell,b_\ell)=\sign(b_\ell)\in\{-1,1\}$  in the strong
field case. Existence and non-vanishing of the degree is the condition of \cite[Theorem 6.1]{Kn}.\\
$\bullet$
The conjugacy property follows, as ${\cal P}^k({\cal P}(x))\in {\cal H}_{\sigma(a)_{k}}$ is equivalent to
${\cal P}^{k+1}(x)$ $\in {\cal H}_{a_{k+1}}$ and to ${\cal P}^k(x)\in {\cal H}_{a_{k}}$\quad $(k\in\bZ)$.
\hfill $\Box$
\begin{remarks} \quad\\[-6mm]
\begin{enumerate}
\item 
The result is quite different from the one in \cite{KSS}. There, for each of 
the {\em continuous} radially symmetric bumps, one had assumed existence of an (unstable)
periodic orbit. This then lead to arbitrarily large winding numbers.
The symbolic dynamics of scattering then gave rise to a {\em semi\,}conjugacy, see  
\cite[Theorem 3.4]{KSS}.
\item 
Similarly, one may consider trapped orbits (with $\lim_{t\to \pm\infty} \|q(t)\|=\infty$
but $\limsup_{t\to \mp\infty} \|q(t)\|<\infty$ and scattering orbits 
(with $\lim_{t\to +\infty} \|q(t)\|=\lim_{t\to -\infty} \|q(t)\|=\infty$). There one uses half-infinite respectively
finite symbol sequences. Prescribing initial and/or final directions of these orbits, one again
obtains symbolic conjugacies. Then, however, one has to exclude directions that are given
by oriented lines meeting two domains $C_\ell$, $C_m$. 
\hfill $\Diamond$
\end{enumerate}
\end{remarks}

%
%


\begin{thebibliography}{99}
\bibitem[Ga]{Ga}
Pierre Gaspard: 
\emph{Chaos, scattering and statistical mechanics.}
Cambridge University Press,  2005
\bibitem[KK]{KK}
Andreas Knauf, Markus Krapf: 
\emph{The non-trapping degree of scattering.} Nonlinearity {\bf 21}, 2023--2041 (2008)
\bibitem[KSS]{KSS}
Andreas Knauf, Frank Schulz, Karl Friedrich Siburg:
\emph{Positive topological entropy for multi-bump magnetic fields.} 
Nonlinearity {\bf 26}, 727--743 (2013)
\bibitem[Kn]{Kn}
Andreas Knauf: \emph{Qualitative Aspects of Classical Potential Scattering.}
Regular and Chaotic Dynamics,  {\bf 4}, No.1, 1--20 (1999) 
\bibitem[Sm]{Sm}
Uzy Smilansky:  \emph{The Classical and Quantum Theory of Chaotic Scattering.}
Lecture Notes. Les Houches, 1989
\end{thebibliography}
\end{document}